\documentclass{article}
\RequirePackage{amsthm,amsmath,amssymb,amsthm,xcolor,mathabx}
\usepackage{graphicx,dsfont,appendix,caption}
\usepackage[numbers]{natbib}

\usepackage[utf8]{inputenc}

\numberwithin{equation}{section}
\theoremstyle{plain}
\newtheorem{thm}{Theorem}[section]

\newtheorem{pro}{Proposition}
\newtheorem*{hyp}{Hypothesis}

\newenvironment{Figure}
  {\par\medskip\noindent\minipage{\linewidth}}
  {\endminipage\par\medskip}

\title{Raking-ratio empirical process with auxiliary information learning}
\author{Mickael Albertus\footnote{mickael.albertus@math.univ-toulouse.fr}}
\date{}

\begin{document}

\maketitle

\begin{abstract}
	The raking-ratio method is a statistical and computational method which adjusts the empirical measure to match the true probability of sets of a finite partition. We study the asymptotic behavior of the raking-ratio empirical process indexed by a class of functions when the auxiliary information is given by estimates. We suppose that these estimates result from the learning of the probability of sets of partitions from another sample larger than the sample of the statistician, as in the case of two-stage sampling surveys. Under some metric entropy hypothesis and conditions on the size of the information source sample, we establish the strong approximation of this process and show in this case that the weak convergence is the same as the classical raking-ratio empirical process. We also give possible statistical applications of these results like the strengthening of the $Z$-test and the chi-square goodness of fit test.\\

\noindent \begin{it}Keywords:\end{it} Uniform central limit theorems, Nonparametric statistics, empirical processes, raking ratio process, auxiliary information, learning.\\

\noindent \begin{it}MSC Classification:\end{it} 62G09, 62G20, 60F17, 60F05.
\end{abstract}

\section{Introduction}

\textbf{Description.} The raking-ratio method is a statistical and computational method aiming to incorporate auxiliary information given by the knowledge of probability of a set of several partitions. The algorithm modifies a sample frequency table in such a way that the marginal totals satisfy the known auxiliary information. At each turn, the method performs a simple cross-multiplication and assigns new weights to individuals belonging to the same set of a partition in order to satisfy the known constraints: it is the "ratio" step of this method. After each modification, the previous constraints are no longer fulfilled in general. Nevertheless, under the conditions that all initial frequencies are strictly positive, if we iteratively cycle the ratio step through a finite number of partitions, the method converges to a frequency table satisfying the expected values -- see~\cite{Sinkh64}. It is the "raking" step of the algorithm. The goal of these operations is therefore to improve the quality of estimators or the power of statistical tests based on the exploitation of the sample frequency table by lowering the quadratic risk when the sample size is large enough. For a numerical example of the raking-ratio method, see Appendix A.1 of~\cite{AlbBert18}. For an example of a simple statistic using the new weights from the raking-ratio method see Appendix~\ref{Appendix0}. The following paragraph summarizes the known results for this method.\medskip

\noindent \textbf{Literature.} The raking-ratio method was suggested by Deming and Stephan and called in a first time "iterative proportions" -- see Section 5 of~\cite{DemSte40}. This algorithm has been initially proposed to adjust the frequency table in the aim to converge it towards the least squares solution. Stephan~\cite{Ste42} then showed that this last statement was wrong and proposed a modification to correct it. Ireland and Kullback~\cite{IreKull68} proved that the raking-ratio method converges to the unique projection of the empirical measure with Kullback-Leibler distance on the set of discrete probability measures verifying all knowing constraints. In some specific cases, estimates for the variance of cell probabilities in the case of a two-way contingency table were established: Brackstone and Rao~\cite{BraRao79} for $ N \leq 4 $, Konijn~\cite{Kon81} or Choudhry and Lee~\cite{ChouLee87}, Bankier~\cite{Bank86} for $ N = 2 $ and Binder and Th{\'e}berge~\cite{BinThe88} for any $ N $. Results of these papers suggest the decrease of variance for the raked estimators of the cells of the table and for a finite number of iterations by providing a complex approximation of the variance of these estimators. Albertus and Berthet~\cite{AlbBert18} defined the empirical measure and process associated to the raking-ratio method and have proved the asymptotic bias cancellation, the asymptotic reduction of variance and so the diminution of the quadratic risk for these process. To prove it, they showed that the raking-ratio empirical process indexed by a class of functions satisfying some metric entropy conditions converges weakly to a specific centered Gaussian process with a lower variance than the usual Brownian bridge. Under general and natural conditions that are recalled below, they proved that the variance decreases by raking among the same cycle of partitions.\medskip

\noindent \textbf{Auxiliary information learning.} The main motivation of this paper is when the statistician does not have the true probability of sets of a given partition but has a source of information which gives him an estimation of this probability more precisely than if he used his own sample. This source can be of different types: preliminary survey of a large sample of individuals, database processing, purchase of additional data at a lower cost, the knowledge of an expert... We suppose in our model that only the estimate of the auxiliary information is transmitted by the source. This hypothesis ensures a fast speed of data acquisition and allows a plurality of sources of information and a diversity of partitions. It is a common situation in statistics since today's technologies like streaming data allow the collection and the transmission of such information in real time. The statistician can use this learned information as auxiliary information which is an estimate of the true one. The raking-ratio method makes it possible to combine shared information of several sources. The main statistical question of this article is whether the statistician can still apply the raking-ratio method by using the estimate of inclusion probabilities rather than the true ones as auxiliary information. We will show that the answer to this question is positive provided that we control the minimum size of the samples of the different sources of auxiliary information.\medskip

\noindent \textbf{Organization.} This paper is organized as follow. Main notation and results are respectively grouped at Section~\ref{Notation_sec} and Section~\ref{Resul_sec}. Some statistical applications are given at Section~\ref{Stats_sec}. We end up by exposing all the proofs at Section~\ref{Proofs_sec}. Appendix~\ref{Appendix0} contains a numerical example of the calculation of a raked mean on a generated sample. At Appendix~\ref{Appendix1} we do the calculation of the asymptotic variance of the raked Gaussian process in a simple case.

\section{Results of the paper}

\subsection{Main notation}\label{Notation_sec}

\textbf{Framework.} Let $ X_1,\dots,X_n, X $ be i.i.d. random variables defined on the same probability space $(\Omega,\mathcal{T},\mathbb{P})$ with same unknown law $ P=\mathbb{P}^{X_{1}} $ on some measurable space $ (\mathcal{X},\mathcal{A}) $. We endow the measurable space $(\mathcal{X},\mathcal{A})$ with $P$.\medskip

\noindent \textbf{Class of functions.} Let $ \mathcal{M} $ denote the set of real valued measurable functions on $ (\mathcal{X},\mathcal{A}) $. We consider a class of functions $ \mathcal{F} \subset \mathcal{M} $ such that $ \sup_{f \in \mathcal{F}} |f| \leq M_\mathcal{F} < +\infty $ for some $ M_\mathcal{F}>0 $ and satisfying the pointwise measurability condition, that is there exists a countable subset $ \mathcal{F}_* \subset \mathcal{F} $ such that for all $ f \in \mathcal{F} $ there exists a sequence $ \{f_m\} \subset \mathcal{F}_* $ with $ f $ as simple limit, that is $ \lim_{m \to +\infty} f_m(x) = f(x) $ for all $ x \in \mathcal{X} $. This condition is often used to ensure the $ P $-measurability of $ \mathcal{F} $ -- see example 2.3.4 of~\cite{VanWell}. For a probability measure $ Q $ on $ (\mathcal{X},\mathcal{A}) $ and $ f,g \in \mathcal{M} $ let $ d_Q^2(f,g) = \int_\mathcal{X} (f-g)^2 dQ $. Let $ N(\mathcal{F},\varepsilon,d_Q) $ be the minimimum number of balls with $ d_Q $-radius $ \varepsilon $ necessary to cover $ \mathcal{F} $ and $ N_{[\ ]}(\mathcal{F},\varepsilon,d_Q) $ be the least number of $ \varepsilon $-brackets necessary to cover $ \mathcal{F} $, that is elements of the form $ [g_-,g_+] = \{f \in \mathcal{F} : g_- \leq f \leq g_+\} $ with  $ d_P(g_-,g_+) <\varepsilon $. We also assume that $ \mathcal{F} $ satisfies one of the two metric entropy conditions (VC) or (BR) discussed below.

\begin{hyp}[VC]
    For $ c_0,\nu_0>0 $, $ \sup_Q N(\mathcal{F},\varepsilon,d_Q) \leq c_0/\varepsilon^{\nu_0} $ where the supremum is taken over all discrete probability measures $ Q $ on $ (\mathcal{X},\mathcal{A}) $.
\end{hyp}
\begin{hyp}[BR]
For $ b_0 > 0$, $ r_0 \in (0,1) $, $ N_{[ \ ]}(\mathcal{F},\varepsilon,d_P) \leq \exp(b_0^2/\varepsilon^{2r_0}) $.
\end{hyp}

\noindent If we add to $ \mathcal{F} $ all elements $ f \mathds{1}_{A_j^{(N)}} $ for every $ N > 0, 1 \leq j \leq m_N $ and $ f \in \mathcal{F} $, $ \mathcal{F} $ still satisfies the same entropy condition but with a new constant $ c_0 $ or $ b_0 $. We denote $ \ell^\infty(\mathcal{F}) $ the set of real-valued functions bounded on $ \mathcal{F} $ endowed with the supremum norm $ ||\cdot||_\mathcal{F} $. In this paper the following notations are used: for all $ f \in \mathcal{F}, A \in \mathcal{A} $ we denote $ P(f) = \mathbb{E}[f(X)] $, $ P(A) = P(\mathds{1}_A) $, $ \mathbb{E}[f|A] = P(f \mathds{1}_{A})/P(A) $, $ \sigma_f^2 = \mathrm{Var}(f(X)) $ and $ \sigma_\mathcal{F}^2 = \sup_{f \in \mathcal{F}}\sigma_f^2 $.\medskip

\noindent \textbf{Empirical measures and processes.} We denote the empirical measure $ \mathbb{P}_n(\mathcal{F}) = \{\mathbb{P}_n(f) : f \in \mathcal{F}\} $ defined by $ \mathbb{P}_n(f) = \frac{1}{n} \sum_{i=1}^n f(X_i) $ and the empirical process $ \alpha_n(\mathcal{F})  = \{\alpha_n(f) : f \in \mathcal{F}\}$ defined by $ \alpha_n(f) = \sqrt{n}(\mathbb{P}_n(f)-P(f)) $. For $ N \in \mathbb{N} $, let $$ \mathcal{A}^{(N)} = \{A_1^{(N)},\dots,A_{m_N}^{(N)}\}\subset \mathcal{A}, $$ be a partition of $ \mathcal{X} $ such that $$ P[\mathcal{A}^{(N)}] = (P(A_1^{(N)}),\dots,P(A_{m_N}^{(N)})) \neq \mathbf{0}. $$ Let $ \mathbb{P}_n^{(N)}(\mathcal{F}) = \{{\mathbb{P}}_n^{(N)}(f) : f \in \mathcal{F}\} $ be the $ N $-th raking-ratio empirical measure defined recursively by $ \mathbb{P}_n^{(0)} = \mathbb{P}_n $ and for all $ f \in \mathcal{F} $,
\begin{align*}
    \mathbb{P}_n^{(N)}(f) = \sum_{j=1}^{m_N} \frac{P(A_j^{(N)})}{\mathbb{P}_n^{(N-1)}(A_j^{(N)})} \mathbb{P}_n^{(N-1)}(f \mathds{1}_{A_j^{(N)}}).
\end{align*} The empirical measure $ \mathbb{P}_n^{(N)}(\mathcal{F}) $ uses the auxiliary information given by $ P[\mathcal{A}^{(N)}] $ to modify $ \alpha_n(\mathcal{F}) $ such that $$ \mathbb{P}_n^{(N)}[\mathcal{A}^{(N)}] = (\mathbb{P}_n^{(N)}(A_1^{(N)}),\dots,\mathbb{P}_n(A_{m_N}^{(N)})) = P[\mathcal{A}^{(N)}] .$$ We denote $ \alpha_n^{(N)}(\mathcal{F}) = \{\alpha_n^{(N)}(f) : f \in \mathcal{F}\} $ the $N$-th raking-ratio empirical process defined for all $ f \in \mathcal{F} $ by \begin{align}
    \label{DefalphN} \alpha_n^{(N)}(f) = \sqrt{n}(\mathbb{P}_n^{(N)}(f)-P(f)).
\end{align} This process satisfies the following property $$ \alpha_n^{(N)}[\mathcal{A}^{(N)}] = (\alpha_n^{(N)}(A_1^{(N)}),\dots,\alpha_n^{(N)}(A_{m_N}^{(N)})) = \mathbf{0}. $$

\noindent \textbf{Gaussian processes.} Under (VC) or (BR), $ \mathcal{F} $ is a Donsker class, that is $ \alpha_n(\mathcal{F}) $ converges weakly in $ \ell^\infty(\mathcal{F}) $ to the $ P $-Brownian bridge $ \mathbb{G}(\mathcal{F}) = \{ \mathbb{G}(f) : f \in \mathcal{F}\} $, the Gaussian process such that $ f \mapsto \mathbb{G}(f) $ is linear and for all $ f,g \in \mathcal{F} $, $$ \mathbb{E}[\mathbb{G}(f)] = 0, \ \mathrm{Cov}(\mathbb{G}(f),\mathbb{G}(g)) = P(fg)-P(f)P(g). $$ For short, we denote $ \mathbb{G}(A) = \mathbb{G}(\mathds{1}_A) $ for any $ A \in \mathcal{A} $. Let $ \mathbb{G}^{(N)}(\mathcal{F}) = \{ \mathbb{G}^{(N)}(f) : f \in \mathcal{F}\} $ be the $ N $-th raking-ratio $ P $-Brownian bridge, that is a centered Gaussian process defined recursively by $ \mathbb{G}^{(0)} = \mathbb{G} $ and for any $ N > 0, f \in \mathcal{F} $, \begin{align}
    \mathbb{G}^{(N)}(f) = \mathbb{G}^{(N-1)}(f) - \sum_{j=1}^{m_N} \mathbb{E}[f|A_j^{(N)}] \mathbb{G}^{(N-1)}(A_j^{(N)}).\label{DefGN}
\end{align}
Albertus and Berthet established the strong approximation and the weak convergence when $ n$ goes to infinity in $ \ell^\infty(\mathcal{F})$ of $ \alpha_n^{(N)}(\mathcal{F}) $ to $ \mathbb{G}^{(N)}(\mathcal{F}) $ for $ N $ fixed -- see Proposition 4 and Theorem 2.1 of~\cite{AlbBert18}. For that they used the strong approximation of the empirical process indexed by a function class satisfying~(VC) or~(BR) -- see Theorem 1 and 2 of~\cite{Berth06}. They gave the exact value of $ \sigma_f^{(N)} = \mathrm{Var}(\mathbb{G}^{(N)}(f)) $ and showed in particular for all $ f \in \mathcal{F} $ and $ N_0 \in \mathbb{N} $ that $ \sigma_f^{(N_0)} \leq \sigma_f^{(0)}=\sigma_f $ and $\sigma_f^{(N_1)} \leq \sigma_f^{(N_0)} $ if $ N_1 \geq 2N_0 $ is such that $ \mathcal{A}^{(N_0-k)}=\mathcal{A}^{(N_1-k)} $ for $ 0\leq k \leq N_0 $ -- see Propositions 7, 8, 9.\medskip

\noindent \textbf{Auxiliary information.} For $ N > 0 $ let $ \mathbb{P}_N'[\mathcal{A}^{(N)}] = (\mathbb{P}_N'(A_1^{(N)}),\dots,\mathbb{P}_N'(A_{m_N}^{(N)})) $ be a random vector with multinomial law, $ n_N $ trials and event probabilities $ P[\mathcal{A}^{(N)}] $. This random vector corresponds to the estimation of the auxiliary information of the $ N $-th auxiliary information source based on a sample of size $n_N = n_N(n) \gg n $ not necessarily independent of $ X_1,\dots,X_n $. We study the asymptotic behavior of the raking-ratio empirical process which uses $ \mathbb{P}_{N}'[\mathcal{A}^{(N)}] $ as auxiliary information instead of $ P[\mathcal{A}^{(N)}] $. By defining the sequence $ \{n_N\} $ we suppose that this information can be estimated by different sources that would not necessarily have the same sample size but still have a sample size larger than $ n $. Let $ \widetilde{\mathbb{P}}_n^{(N)}(\mathcal{F}) = \{ \widetilde{\mathbb{P}}_n^{(N)}(f) : f \in \mathcal{F}\} $ be the $ N$-th raking-ratio empirical measure with learned auxiliary information defined recursively by $ \widetilde{\mathbb{P}}_n^{(0)} = \mathbb{P}_n $ and for all $ N > 0, f \in \mathcal{F} $,
\begin{align*}
    \widetilde{\mathbb{P}}_n^{(N)}(f) = \sum_{j=1}^{m_N} \frac{\mathbb{P}_{N}'(A_j^{(N)})}{\widetilde{\mathbb{P}}_n^{(N-1)}(A_j^{(N)})} \widetilde{\mathbb{P}}_n^{(N-1)}(f \mathds{1}_{A_j^{(N)}}).
\end{align*} This empirical measure satisfies the learned auxiliary information since \begin{align*}
    \widetilde{\mathbb{P}}_n^{(N)}[\mathcal{A}^{(N)}] &= (\widetilde{\mathbb{P}}_n^{(N)}(A_1^{(N)}),\dots,\widetilde{\mathbb{P}}_n^{(N)}(A_{m_N}^{(N)}) ) \\
    &= \mathbb{P}_N'[\mathcal{A}^{(N)}].
\end{align*} We define $ \widetilde{\alpha}_n^{(N)}(\mathcal{F}) = \{\widetilde{\alpha}_n^{(N)}(f) : f \in \mathcal{F} \} $ the $ N $-th raking-ratio empirical with estimated auxiliary information defined for $ f \in \mathcal{F} $ by \begin{align}
    \label{DefalphtN}  \widetilde{\alpha}_n^{(N)}(f) = \sqrt{n}(\widetilde{\mathbb{P}}_n^{(N)}(f)-P(f)).
\end{align}

\subsection{Main results}\label{Resul_sec}
 For $ N_0 > 0 $, denote $ K_\mathcal{F} = \max(1,M_\mathcal{F}) $ and \begin{align*}
     p_{(N_0)} &= \min_{1 \leq N \leq N_0} \min_{1 \leq j \leq m_N} P(A_j^{(N)}), \\
     m_{(N_0)} &= \sup_{0 \leq N \leq N_0} m_N, \\
     n_{(N_0)} &= \min_{1 \leq N \leq N_0} n_N > n.
 \end{align*}  Empirical measures $ \mathbb{P}_n^{(0)}(\mathcal{F}),\dots,\mathbb{P}_n^{(N_0)}(\mathcal{F}) $ and $ \widetilde{\mathbb{P}}_n^{(0)}(\mathcal{F}),\dots,\widetilde{\mathbb{P}}_n^{(N)}(\mathcal{F}) $ are defined on the set $$ B_{n,N_0} = \left\{ \min_{0 \leq N \leq N_0} \min_{1\leq j \leq m_N} \mathbb{P}_n(A_j^{(N)}) > 0 \right\}, $$ which satisfies $$ \mathbb{P}(B_{n,N_0}^C) \leq \sum_{N = 1}^{N_0} m_N (1-p_N)^n \leq N_0 m_{(N_0)} (1-p_{(N_0)})^n,  $$ where $ B_{n,N_0}^C = \Omega \setminus B_{n,N_0}. $
 The following proposition bounds the probability that $ ||\widetilde{\alpha}_n^{(N)}||_\mathcal{F} $ deviates from a certain value.
\begin{pro}\label{ProTalIneq}
    For any $ N_0 \in \mathbb{N} $, $ n > 0 $ and $ t>0 $, it holds under the event $ B_{n,N_0}$ \begin{align}
\mathbb{P}\left(\sup_{0 \leq N \leq N_0} ||\widetilde{\alpha}_n^{(N)}||_\mathcal{F} > t\right) &\leq N_0 \mathbb{P}\left( ||\widetilde{\alpha}_n^{(0)}||_\mathcal{F} > \frac{t p_{(N_0)}^{N_0}}{4^{N_0} m_{(N_0)}^{N_0} K_\mathcal{F}^{N_0} (1+t/\sqrt{n})^{N_0}} \right) \notag \\
&\quad + 2N_0^3 m_{(N_0)}\exp\!\left(-\frac{ n_{(N_0)} p_{(N_0)}^2 t^2}{2n m_{(N_0)}^2 K_\mathcal{F}^2} \right).\label{BoundalphtN_Result}
\end{align}
Under (VC) and the event $ B_{n,N_0}$ there exists $ t_0 >0 $ such that for all $ t_0 < t <2M_\mathcal{F}\sqrt{n} $,  \begin{align}
    \mathbb{P}\left(\sup_{0 \leq N \leq N_0} ||\widetilde{\alpha}_n^{(N)}||_\mathcal{F} > t\right) &\leq D_1 t^{\nu_0}\exp(-D_2 t^2) \notag \\
    &\quad+ 2N_0^3 m_{(N_0)}\exp\!\left(-\frac{n_{(N_0)} p_{(N_0)}^2 t^2}{2n m_{(N_0)}^2 K_\mathcal{F}^2} \right),  \label{BoundalphtN_ResultVC}
\end{align} where $ D_1,D_2>0 $ are defined by~\eqref{DefD1D2}. Under (BR) and the event $ B_{n,N_0}$ there exists $ t_0,C > 0 $ such that for all $ t_0<t<C\sqrt{n} $, \begin{align}
    \mathbb{P}\left( \sup_{0 \leq N \leq N_0} ||\widetilde{\alpha}_n^{(N)}||_\mathcal{F} > t\right) &\leq D_3 \exp(-D_4 t^2) \notag\\
    &\quad+ 2N_0^3 m_{(N_0)}\exp\!\left(-\frac{n_{(N_0)} p_{(N_0)}^2 t^2}{2n m_{(N_0)}^2 K_\mathcal{F}^2} \right),\label{BoundalphtN_ResultBR}
\end{align} where $ D_3,D_4 >0 $ are defined by~\eqref{DefD3D4}.
\end{pro} \noindent Proposition~\ref{ProTalIneq} proves that if $ \mathcal{F} $ satisfies~(VC) or~(BR) then almost surely $ ||\alpha_n||_\mathcal{F} =O(\sqrt{\log(n)}) $.  If $ \mathcal{F} $ satisfies (VC), let define $ v_n = n^{-\alpha_0}(\log n)^{\beta_0} $ with $ \alpha_0=1/(2+5\nu_0) \in (0,1/2) $ and $ \beta_0=(4+5\nu_0)/(4+10\nu_0) $. If $ \mathcal{F} $ satisfies (BR), let define $ v_n = (\log n)^{-\gamma_0} $ with $ \gamma_0=(1-r_0)/2r_0 $. The following result establishes the strong approximation of  $ \widetilde{\alpha}_n^{(N)}(\mathcal{F}) $ by $ \mathbb{G}^{(N)}(\mathcal{F}) $.
\begin{thm}\label{ThmApproxForte}
     Let $ N_0 \in \mathbb{N} $. There exists $ d_0,n_0>0 $, a sequence $ \{X_n\} $ of independent random variables with law $ P $ and a sequence $ \{\mathbb{G}_n\} $ of versions of $ \mathbb{G} $ supported on a same probability space such that for all $ n > n_0 $,
       \begin{align}
        \mathbb{P}\left(\sup_{0 \leq N \leq N_0} ||\widetilde{\alpha}_n^{(N)}-\mathbb{G}_n^{(N)}||_\mathcal{F} > d_0 \! \left(v_n\!+\!\sqrt{\frac{n\log(n)}{n_{(N_0)}}} \right)\!\right)\!< \!\frac{1}{n^2}, \label{ThmApproxForteResult}
    \end{align} where $ \mathbb{G}_n^{(N)} $ is the version of $ \mathbb{G}^{(N)}$ derived from $ \mathbb{G}_n^{(0)} = \mathbb{G}_n $ through~\eqref{DefGN}.
\end{thm}
\noindent By Borel-Cantelli lemma we have almost surely for large $ n $, \begin{align}\label{ThmApproxForteResultBis}
    \sup_{0 \leq N \leq N_0} ||\widetilde{\alpha}_n^{(N)}-\mathbb{G}_n^{(N)}||_\mathcal{F} \leq d_2 \! \left(v_n + \sqrt{\frac{n\log(n)}{n_{(N_0)}}}\right).
\end{align} Sequence $ v_n $ in the previous bound is the deviation from $ \alpha_n^{(N)}(\mathcal{F}) $ to $ \mathbb{G}_n^{(N)}(\mathcal{F}) $ while $ \displaystyle \sqrt{n\log(n)/n_{(N_0)}} $ represents the deviation from $ \widetilde{\alpha}_n^{(N)}(\mathcal{F}) $ to $ \alpha_n^{(N)}(\mathcal{F}) $. Under the condition that the sample size of the sources are large enough, Theorem~\ref{ThmApproxForte} implies that the sequence $ (\widetilde{\alpha}_n^{(0)}(\mathcal{F}),\dots,\widetilde{\alpha}_n^{(N_0)}(\mathcal{F})) $ converges weakly to $ (\mathbb{G}^{(0)}(\mathcal{F}),\dots,\mathbb{G}^{(N_0)}(\mathcal{F})) $ on $ \ell^\infty(\mathcal{F} \to \mathbb{R}^{N_0+1}) $ as the same way as $ (\alpha_n^{(0)}(\mathcal{F}),\dots,\alpha_n^{(N_0)}(\mathcal{F})) $.

\subsection{Statistical applications}\label{Stats_sec}

\textbf{Improvement of a statistical test.} Any statistical test using the empirical process can be adapted to use auxiliary information to strengthen this test. It suffices to replace in the expression of the test statistic the process $ \alpha_n(\mathcal{F}) $ by $ \alpha_n^{(N)}(\mathcal{F}) $ if we have the true auxiliary information or by $ \widetilde{\alpha}_n^{(N)}(\mathcal{F}) $ if we have an estimation of this information. The two following subsections give an example of application in the case of the $ Z $-test and the chi-squared goodness of fit test. In both case, we transform the statistic of theses tests and keep the same decision procedure. In the first case, we show that this new statistical test has the same significance level but a higher power. For the second case, we prove that the confidence level decreases and that under $ (H_1) $, the new statistic goes to infinity as the same way as the usual one.\medskip

\noindent \textbf{$Z$-test.} This test is used to compare the mean of a sample to a given value when the variance of the sample is known. The null hypothesis is $ (H_0): P(f) = P_0(f), $ for some $ f \in \mathcal{F} $ and a probability measure $ P_0 \in \ell^\infty(\mathcal{F}) $. The statistic of the classical $Z$-test is $$ Z_n = \sqrt{n}\frac{\mathbb{P}_n(f) - P_0(f)}{\sigma_f}. $$ Under $ (H_0) $, asymptotically the statistic $ Z_n $ follows the standard normal distribution. We reject the null hypothesis at the $ \alpha $ level when $ |Z_n| > t_\alpha $, $ t_\alpha = \Phi(1-\alpha/2)  $ with $ \Phi $ the probit function. Let define the following statistics \begin{align*}
    Z_n^{(N)} &= \sqrt{n}\frac{\mathbb{P}_n^{(N)}(f) - P_0(f)}{\sigma_f^{(N)}}, \\
    \widetilde{Z}_n^{(N)} &= \sqrt{n}\frac{\widetilde{\mathbb{P}}_n^{(N)}(f) - P_0(f)}{\sigma_f^{(N)}},
\end{align*} Since the law $ P $ is unknown, $ \sigma_f $ and $ \sigma_f^{(N)}$ for $ N \geq 1 $ are usually unknown but a consistent estimation of these variances can be used to calculate $ Z_n, Z_n^{(N)} $ or $ \widetilde{Z}_n^{(N)} $ -- a concrete example of this remark is given at the following paragraph. Doing it does not change the asymptotic behavior of the random variables $ Z_n, Z_n^{(N)} $ and $ \widetilde{Z}_n^{(N)} $, whether the hypothesis $ (H_0) $ is verified or not. The statistical tests based on the reject decision $ |Z_n^{(N)}| > t_\alpha $ and $ |\widetilde{Z}_n^{(N)}|>t_\alpha $ have the same significance level than the usual test based on the decision $ |Z_n|>t_\alpha $ since, under $ (H_0)$, $ Z_n^{(N)} $ and $ \widetilde{Z}_n^{(N)} $ converge weakly to $ \mathcal{N}(0,1)$ -- see Proposition 6 of~\cite{AlbBert18}. The following proposition shows that the ratio of the beta risk of the usual $ Z$-test and the new statistical test with auxiliary information goes to infinity as $ n \to +\infty $. \begin{pro}\label{pro_Ztest}
    Assume that $\sigma_f^{(N)} < \sigma_f $. Under $ (H_1) $, for all $ \alpha \in (0,1) $ and $ n$ large enough one have \begin{align}
        \frac{\mathbb{P}(|Z_n| \leq t_\alpha)}{\mathbb{P}(|Z_n^{(N)}| \leq t_\alpha)} \geq \exp \left(n(P(f)-P_0(f))^2 \left(\frac{1}{\sigma_f^{(N)}}-\frac{1}{\sigma_f} \right) \right).\label{pro_Ztest_3}
    \end{align}
\end{pro}

\noindent \textbf{Z-test in a simple case.} To calculate $ Z_n^{(N)} $ or $ \widetilde{Z}_n^{(N)} $ one needs the expression of $ \sigma_f^{(N)} $. To illustrate how to get it we work on a simple case, when the auxiliary information is given by probabilities of two partitions of two sets. More formally for $ k \in \mathbb{N}^* $ we define $ \mathcal{A}^{(2k-1)} = \mathcal{A}=\{A,A^C\} $ and $ \mathcal{A}^{(2k)} =\mathcal{B}= \{B,B^C\} $. By using Proposition 7 of~\cite{AlbBert18} we give simple expressions of $ \sigma_f^{(N)} $ for $ N=1,2 $.  For the sake of simplification, let denote \begin{align}
    p_A &= P(A), \quad p_{\overline{A}} = P(A^C),\quad p_B=P(B), \quad p_{\overline{B}} = P(B^C),\notag \\
    p_{AB}&=P(A\cap B), \quad \Delta_A = \mathbb{E}[f|A]-\mathbb{E}[f], \quad  \Delta_B = \mathbb{E}[f|B]-\mathbb{E}[f], \label{ShortcutVar}
\end{align} then, \begin{align*}
    \sigma_f^{(1)} &= \sigma_f - \mathbb{E}[f|\mathcal{A}]^t \cdot \mathrm{Var}(\mathbb{G}[\mathcal{A}]) \cdot \mathbb{E}[f|\mathcal{A}] \\
    &= \sigma_f - p_A p_{\overline{A}} (\mathbb{E}[f|A]-\mathbb{E}[f|A^C])^2, \\
    \sigma_f^{(2)} &= \sigma_f- \mathbb{E}[f|\mathcal{B}]^t \cdot \mathrm{Var}(\mathbb{G}[\mathcal{B}]) \cdot \mathbb{E}[f|\mathcal{B}] \\
    &\quad - \left(\mathbb{E}[f|\mathcal{A}] - \mathbf{P}_{\mathcal{B}|\mathcal{A}} \cdot \mathbb{E}[f|\mathcal{B}] \right)^t \cdot \mathrm{Var}(\mathbb{G}[\mathcal{A}]) \cdot  \left(\mathbb{E}[f|\mathcal{A}]- \mathbf{P}_{\mathcal{B}|\mathcal{A}} \cdot \mathbb{E}[f|\mathcal{B}] \right) \\
    &=\sigma_f -   p_B p_{\overline{B}} (\mathbb{E}[f|B]-\mathbb{E}[f|B^C])^2\\
    &\quad - \left( p_A p_{\overline{A}} +\frac{p_B p_{\overline{B}} (p_{AB}-p_A p_B)}{p_A^2 p_{\overline{A}}^2} \right) (\mathbb{E}[f|A]-\mathbb{E}[f|A^C])^2 ,
\end{align*} where $ \mathbf{P}_{\mathcal{A}|\mathcal{B}},\mathbf{P}_{\mathcal{B}|\mathcal{A}} $ are stochastic matrices given by~\eqref{DefPABPBA}, $ \mathbb{E}[f|\mathcal{A}],\mathbb{E}[f|\mathcal{B}] $ are conditional expectation vectors given by~\eqref{DefEAEB} and $ \mathrm{Var}(\mathbb{G}[\mathcal{A}]), \mathrm{Var}(\mathbb{G}[\mathcal{B}]) $ are the covariance matrices of $ \mathbb{G}[\mathcal{A}]=(\mathbb{G}(A),\mathbb{G}(A^C)) $ and $\mathbb{G}[\mathcal{B}] = (\mathbb{G}(B),\mathbb{G}(B^C)) $ that is the matrices given by~\eqref{DefCovGAGB}. Albertus and Berthet proved that the raked Gaussian process $ \mathbb{G}^{(N)} $ converges almost surely as $ N \to +\infty $ to some centered Gaussian process $ \mathbb{G}^{(\infty)} $ with an explicit expression. The stabilization of the raking-ratio method in the case of two marginals when $ N \to +\infty $ is fast since the Levy-Prokhorov distance between $ \mathbb{G}^{(N)} $ and $ \mathbb{G}^{(\infty)} $ is almost surely at most $ O(N \lambda^{N/2}) $ for some $ \lambda \in (0,1) $ -- see Proposition 11 of~\cite{AlbBert18}. We denote $ \mathbb{P}_n^{(\infty)}(\mathcal{F}) $ the raked empirical measure after stabilization of the raking-ratio algorithm and $ \sigma_f^{(\infty)} = \mathrm{Var}(\mathbb{G}^{(\infty)}(f)) $ the asymptotic variance. Let define the following statistic $$ Z_n^{(\infty)} = \sqrt{n} \frac{\mathbb{P}_n^{(\infty)}(f)-P_0(f)}{\sigma_f^{(\infty)}}. $$ According to Proposition~\ref{pro_Ztest}, the statistical test based on the reject decision $ |Z_n^{(\infty)}|>t_\alpha $ has the same significance level than the usual $Z$-test based on $ |Z_n|>t_\alpha $ but it is more powerful as $ n $ goes to infinity. In the case of two marginals with two partitions, one can give an explicit and simple expression of the asymptotic variance. By using the notations of~\eqref{ShortcutVar} one have \begin{align}\label{DefSigInf}
    \sigma_f^{(\infty)} =\sigma_f^2- \frac{p_A p_B \left( p_A \Delta_A^2+p_B \Delta_B^2-p_A p_B (\Delta_A-\Delta_B)^2-2p_{AB} \Delta_A \Delta_B \right)}{p_A p_B p_{\overline{A}} p_{\overline{B}}-(p_{AB}-p_A p_B)^2}.
\end{align} The calculation of this variance needs the expression of  $ \mathbb{G}^{(\infty)} $ so it is made at Appendix~\ref{Appendix1}. If we do not have the values given by~\eqref{ShortcutVar} one can use their consistent estimators to estimate the value of $ \sigma_f^{(\infty)} $. If $ \Delta_A=\Delta_B=0 $ then naturally the auxiliary information is useless since $ \sigma_f^{(\infty)} = \sigma_f $, so there is no reduction of the quadratic risk. If $ A $ is independent of $ B $ then $ p_{AB} = p_A p_B $ and $$ \sigma_f^{(\infty)} = \sigma_f- \left( \frac{p_A}{p_{\overline{A}}} \Delta_A^2 + \frac{p_B}{p_{\overline{B}}} \Delta_B^2 \right). $$ \medskip

\noindent \textbf{Chi-square test.} The chi-squared goodness of fit test consists of knowing whether the sample data corresponds to a hypothesized distribution when we have one categorical variable. Let $ \mathcal{B}=\{B_1,\dots,B_m\} $ be a partition of $ \mathcal{X} $. The null hypothesis is \begin{align}
    (H_0): P[\mathcal{B}]=P_0[\mathcal{B}], \label{null_hyp_khi2}
\end{align} where $ P[\mathcal{B}] = (P(B_1),\dots,P(B_m)) $ and $ P_0[\mathcal{B}] = (P_0(B_1),\dots,P_0(B_m)) $, for some probability measure $ P_0 $. The statistic of the classical chi-squared test is $$ T_n = n\sum_{i=1}^m \frac{\left(\mathbb{P}_n(B_i)-P_0(B_i)\right)^2}{P_0(B_i)}. $$ Under $ (H_0) $, asymptotically the statistic $T_n $ follows the $ \chi^2 $ distribution with $ m-1 $ degrees of freedom. We reject the null hypothesis at the level $ \alpha $ when $ Z_n > t_\alpha^{(m)} $, $ t_\alpha^{(m)} = \Phi_{m}(1-\alpha) $ where $ \Phi_m $ is the quantile function of $ \chi^2(m) $. We want to know if the following statistics \begin{align*}
    T_n^{(N)} = n\sum_{i=1}^m \frac{(\mathbb{P}_n^{(N)}(B_i)-P_0(B_i))^2}{P_0(B_i)}, \\
    \widetilde{T}_n^{(N)} =n \sum_{i=1}^m \frac{(\widetilde{\mathbb{P}}_n^{(N)}(B_i)-P_0(B_i))^2}{P_0(B_i)},
\end{align*} somehow improve the test. The following proposition shows that the power of the test is improved with these new statistics. \begin{pro}\label{pro_chi2}
    Under $ (H_0) $ and for all $ \alpha > 0 $, \begin{align}
        \lim_{n\to+\infty} \mathbb{P}(T_n^{(N)}>t_\alpha^{(m)}) &\leq\lim_{n \to +\infty} \mathbb{P}(T_n>t_\alpha^{(m)}) =\alpha, \label{pro_chi2_1}
    \end{align} and if $ n\log(n) = o(n_{(N)}) $ then \begin{align}
        \lim_{n\to+\infty} \mathbb{P}(\widetilde{T}_n^{(N)}>t_\alpha^{(m)}) &\leq \alpha. \label{pro_chi2_2}
    \end{align} Under $ (H_1) $ and for all $ \alpha > 0 $, almost surely there exists $ n_0>0 $ such that for all $ n > n_0 $, \begin{align}
        \min(|T_n|, |T_n^{(N)}|, |T_n^{(N)}|) > t_\alpha^{(m)}.\label{pro_chi2_3}
    \end{align}
\end{pro} \noindent Figure~\ref{fig-loi-khi2} is a numerical example of Proposition~\ref{pro_chi2} under $ (H_0) $. We simulate a two-way contingency table with fixed probabilities $ P[\mathcal{B}], P[\mathcal{A}] $ and we apply the chi-square test with the null hypothesis~\eqref{null_hyp_khi2}. With Monte-Carlo method, we simulate the law of $ T_n $ for $ n=1000$ and the law of $ T_n^{(1)} $ with the auxiliary information given by $ P[\mathcal{A}] $.

\begin{Figure}
\vskip 0.2in
\begin{center}
\includegraphics[scale=0.5]{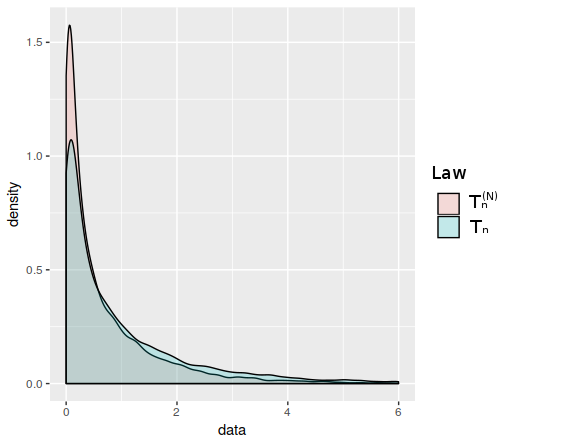}
\label{fig-loi-khi2}
\end{center}
\vskip -0.2in
\end{Figure}\medskip

\noindent \textbf{Costing data.} Another possible statistical application is to study how to share resources -- economic resource, temporal resource, material resource, ... -- to learn auxiliary information from inexpensive data in order to improve the study of statistics on expensive objects. More formally we have a budget $ B $, for our estimates we can buy an individual $ X_i $ at a fixed price $C>0 $ and for the estimation of auxiliary information $ P[\mathcal{A}^{(N)}], N=1,\dots,N_0 $, we can buy the information $ \mathbb{P}_N'[\mathcal{A}^{(N)}] $ at a price $ c_N n_N $ where $ c_N $ is the price for one individual far less than $ C $. The objective is therefore to minimize the bound $ v_n+\sqrt{n\log(n)/n_{(N_0)}} $ proposed by Theorem~\ref{ThmApproxForte} by choosing $n$ high-cost individuals and the $n_1,\dots,n_{N_0} $ low-cost individuals while respecting the imposed budget. So we have to satisfy the following constraint \begin{align}
    Cn+c_1 n_1+\dots+c_{N_0} n_{N_0} \leq B. \label{cost_data_1}
\end{align}
To simplify the problem we will suppose that for all $ 1 \leq N \leq N_0 $, $n_N = n_0 $ and $ c_N = c_0/N_0 $ for some $ c_0 > 0$. It is the case if one pay the auxiliary information from the same auxiliary information source and if one pay all $N_0$ information only once time. Inequality~\eqref{cost_data_1} becomes \begin{align}
    Cn+c_0 n_0 \leq B.\label{cost_data_2}
\end{align} There are several ways to answer this problem. If we want only the strong approximation rate of $ \alpha_n^{(N)} $ by $ \mathbb{G}^{(N)} $ dominates in the uniform error of~\eqref{ThmApproxForteResultBis}, we have to choose $n_0 $ such that $ n_0 \geq n\log(n)/v_n^2 $. If we take $ n_0 = \lceil n\log(n)/v_n^2\rceil $ we could find the maximum value of $ n $ satisfying~\eqref{cost_data_2}. Since $ v_n > \sqrt{\log(n)/n} $ we know that \begin{align}
    n \geq n_{\min}=\left\lfloor  \frac{\sqrt{C^2+4c_0 B} -C}{2c_0}\right\rfloor. \label{n_opt}
\end{align} If we have no way of finding the optimal $n$ -- if we do not have the rate $v_n$ or if we want to avoid additional calculations -- we can take $ n=n_{\min} $ and $ n_0 = \lfloor (B-Cn)/c_0 \rfloor  $ if one want to use the entire budget or $ n_0 =\lceil n\log(n)/v_n^2\rceil $ otherwise.


%



%


\section{Proof} \label{Proofs_sec}

For all this section let fix $ N_0 > 0 $ and let $ \Lambda_n,\Lambda_n'> 0 $ be the following supremum deviations \begin{align*}
    \Lambda_n &= \max\left(\sup_{0 \leq N \leq N_0} ||\widetilde{\alpha}_n^{(N)}||_\mathcal{F}, \sup_{0 \leq N \leq N_0} ||\alpha_n^{(N)}||_\mathcal{F}\right), \\
    \Lambda_n' &= \sup_{1 \leq N \leq N_0} \sup_{1 \leq j \leq m_N} |\alpha_{N}'(A_j^{(N)})|,
\end{align*} where $ \alpha_N'(A_j^{(N)}) = \sqrt{n_N}(\mathbb{P}_N'(A_j^{(N)})-P(A_j^{(N)})) $. Immediately, by Hoeffding inequality we have for all $ \lambda > 0 $,
\begin{align}
	\mathbb{P}\left(\Lambda_n' > \lambda \right) \leq  2N_0 m_{(N_0)} \exp\!\left( -2\lambda^2 \right).\label{MajHoeff}
\end{align}  Now, we give useful decomposition of $ \alpha_n^{(N)}(\mathcal{F}) $ and $ \widetilde{\alpha}_n^{(N)}(\mathcal{F}) $ which will be used in the following proofs. By using definition~\eqref{DefalphN} of $ \alpha_n^{(N)}(\mathcal{F}) $ we have

 \begin{align}
     \alpha_n^{(N)}(f) &= \sqrt{n}\left( \sum_{j=1}^{m_N} \frac{P(A_j^{(N)})}{\mathbb{P}_n^{(N-1)}(A_j^{(N)})} \mathbb{P}_n^{(N-1)}(f \mathds{1}_{A_j^{(N)}}) - P(f \mathds{1}_{A_j^{(N)}}) \right) \notag \\
     &= \sum_{j=1}^{m_N} \frac{P(A_j^{(N)})\alpha_n^{(N-1)}(f\mathds{1}_{A_j^{(N)}})-P(f \mathds{1}_{A_j^{(N)}})\alpha_n^{(N-1)}(A_j^{(N)})}{\mathbb{P}_n^{(N-1)}(A_j^{(N)})}. \label{DecompalphN}
 \end{align}
  As the same way, by using~\eqref{DefalphtN} we have  \begin{align}
     \widetilde{\alpha}_n^{(N)}(f)&= \sum_{j=1}^{m_N} \frac{\mathbb{P}_{N}'(A_j^{(N)})}{\widetilde{\mathbb{P}}_n^{(N-1)}(A_j^{(N)})} \widetilde{\alpha}_n^{(N-1)}(f \mathds{1}_{A_j^{(N)}}) \notag \\
     &- \frac{P(f\mathds{1}_{A_j^{(N)}})}{\widetilde{\mathbb{P}}_n^{(N-1)}(A_j^{(N)})} \left(\widetilde{\alpha}_n^{(N-1)}(A_j^{(N)}) - \sqrt{\frac{n}{n_N}} \alpha_{N}'(A_j^{(N)}) \right).  \label{DecompalphtN}
 \end{align}

\subsection{Proof of Proposition~\ref{ProTalIneq}}

We prove~\eqref{BoundalphtN_Result},~\eqref{BoundalphtN_ResultVC} and~\eqref{BoundalphtN_ResultBR} respectively at Step 1, Step 2 and Step 3.\medskip

\noindent \textbf{Step 1.} Let $ 0 \leq N \leq N_0 $. With~\eqref{DecompalphtN} one can write that
 \begin{align}
	&\mathbb{P}(||\widetilde{\alpha}_n^{(N)}||_\mathcal{F} > t) \notag \\
	&\leq \mathbb{P}\left( \frac{K_\mathcal{F} m_{(N)}\left(2||\widetilde{\alpha}_n^{(N-1)}||_\mathcal{F} + \sqrt{\frac{n}{n_{(N)}}} \Lambda_n'\right)}{p_{(N)} -||\widetilde{\alpha}_n^{(N-1)}||_\mathcal{F}/\sqrt{n}}  > t \right) \notag \\
	&\leq \mathbb{P}\left(\Lambda_n' > \sqrt{\frac{n_{(N)}}{n}} \frac{ tp_{(N)}}{2m_{(N)} K_\mathcal{F}} \right) \notag \\
	&\quad+ \mathbb{P}\left( ||\widetilde{\alpha}_n^{(N-1)}||_\mathcal{F} > \frac{tp_{(N)}}{4m_{(N)}K_\mathcal{F}(1+t/\sqrt{n}) } \right)\notag \\
	&\leq \mathbb{P}\left(\Lambda_n' > \sqrt{\frac{n_{(N_0)}}{n}} \frac{ tp_{(N_0)}}{2m_{(N_0)} K_\mathcal{F}} \right) \notag \\
	&\quad+ \mathbb{P}\left( ||\widetilde{\alpha}_n^{(N-1)}||_\mathcal{F} > \frac{tp_{(N)}}{4m_{(N)}K_\mathcal{F}(1+t/\sqrt{n}) } \right).\label{DecompIneqalphtN}
\end{align}
 By~\eqref{MajHoeff} and induction on~\eqref{DecompIneqalphtN}, we find \begin{align*}
\mathbb{P}\left(||\widetilde{\alpha}_n^{(N)}||_\mathcal{F} > t\right) &\leq \mathbb{P}\left( ||\widetilde{\alpha}_n^{(0)}||_\mathcal{F} > \frac{t p_{(N)}^N}{4^N m_{(N)}^N K_\mathcal{F}^N (1+t/\sqrt{n})^N} \right)\\
&\quad + 2N_0^2 m_{(N_0)}\exp\!\left(-\frac{n_{(N_0)} p_{(N_0)}^2 t^2}{2n m_{(N_0)}^2 K_\mathcal{F}^2} \right).
\end{align*}
The right-hand side of the last inequality is increasing with $ N $ which leads to~\eqref{BoundalphtN_Result}. Since \begin{align}
    \widetilde{\alpha}_n^{(0)}(\mathcal{F}) = \alpha_n(\mathcal{F}) = \alpha_n^{(0)}(\mathcal{F}), \label{EqProcN=0}
\end{align} we can apply Talagrand inequality to control the deviation probability of $ ||\widetilde{\alpha}_n^{(0)}||_\mathcal{F} $ as described in the next two steps.\medskip

\noindent \textbf{Step 2.} According to Theorem 2.14.25 of~\cite{VanWell}, if $ \mathcal{F} $ satisfies (VC) there exists a constant $ D=D(c_0) >0 $ such that, for $ t_0 $ large enough and $ t \geq t_0 $, \begin{align}
    \label{BoundalphaVC}\mathbb{P}\left( ||\widetilde{\alpha}_n^{(0)}||_\mathcal{F} > t \right) \leq \left( \frac{Dt}{M_\mathcal{F} \sqrt{\nu_0}} \right)^{\nu_0} \exp\!\left( \frac{-2t^2}{M_\mathcal{F}^2} \right).
\end{align} Inequalities~\eqref{BoundalphtN_Result} and~\eqref{BoundalphaVC} imply~\eqref{BoundalphtN_ResultVC} for all $ t_0 \leq t \leq 2M_\mathcal{F}\sqrt{n} $,  where $D_1,D_2>0 $ are defined by \begin{align}
    D_1 &= N_0\left( \frac{Dp_{(N_0)}^{N_0}}{\nu_0 4^{N_0} m_{(N_0)}^{N_0} K_\mathcal{F}^{N_0+1}} \right)^{\nu_0}, \notag \\
    D_2 &= \frac{p_{(N_0)}^{2N_0}}{72^{N_0} m_{(N_0)}^{2N_0} K_\mathcal{F}^{3N_0+1}}.\label{DefD1D2}
\end{align}

\noindent \textbf{Step 3.} According to Theorems 2.14.2 and 2.14.25 of~\cite{VanWell}, if $ \mathcal{F} $ satisfies (BR), there exists universal constants $ D,D' > 0 $ such that for all $ t_0  < t < t_1  $, \begin{align}
    \mathbb{P}\left( ||\widetilde{\alpha}_n^{(0)}||_\mathcal{F} > t \right) \leq \exp(-D'' t^2), \label{BoundalphaBR}
\end{align} where $ t_0 =2DM_\mathcal{F}(1+b_0/(1-r_0)) $, $ t_1 = 2D\sigma_\mathcal{F}^2 \sqrt{n}/M_\mathcal{F}$, $ D'' = D'/4D^2 \sigma_\mathcal{F}^2 $. Therefore~\eqref{BoundalphtN_Result} and~\eqref{BoundalphaBR} yields~\eqref{BoundalphtN_ResultBR} where $ D_3,D_4 > 0 $ are defined by \begin{align}
    D_3 = N_0, \quad D_4 = \frac{D'' p_{(N_0)}^{2N_0}}{8^{N_0} m_{(N_0)}^{2N_0} K_\mathcal{F}^{2N_0} (1+2D\sigma_\mathcal{F}^2/M_\mathcal{F})^{2N_0}}. \label{DefD3D4}
\end{align}

\subsection{Proof of Theorem~\ref{ThmApproxForte}}
According to Proposition~\ref{ProTalIneq}, inequality~\eqref{MajHoeff} and Proposition 3 of~\cite{AlbBert18}, there exists $ D > 0 $ such that \begin{align}
    \mathbb{P}\left(\{ \Lambda_n > D \sqrt{\log(n)} \} \ \bigcup \ \{ \Lambda_n' > D\sqrt{\log(n)} \} \right) \leq \frac{1}{3n^2}. \label{BoundDevalphaNAndalphatN}
\end{align}
 According to Theorem 2.1 of~\cite{AlbBert18}, one can define on the same probability space a sequence $ \{X_n\} $ of independent random variable with law $ P $ and a sequence $ \{\mathbb{G}_n\} $ of versions of $ \mathbb{G} $ satisfying the following property. There exists $ n_1, d_1 > 0 $ such that for all $ n > n_1 $, \begin{align*}
    \mathbb{P}\left(\sup_{0 \leq N \leq N_0} ||\alpha_n^{(N)}-\mathbb{G}_n^{(N)}||_\mathcal{F} > d_1 v_n \right) \leq \frac{1}{3n^2},
\end{align*} where $ \mathbb{G}_n^{(N)} $ is the version of $ \mathbb{G}^{(N)}$ derived from $ \mathbb{G}_n^{(0)}=\mathbb{G}_n $ through~\eqref{DefGN}. To show~\eqref{ThmApproxForteResult} it remains to prove, by~\eqref{EqProcN=0}, that for all $ n $ large enough and some $ d_0 > 0 $, $$ \mathbb{P}\left(\sup_{0 \leq N \leq N_0} ||\widetilde{\alpha}_n^{(N)}-\alpha_n^{(N)}||_\mathcal{F} >  d_0\sqrt{\frac{n \log(n)}{n_{(N_0)}}} \right) \leq \frac{2}{3n^2} .$$
Let $ 1 \leq N \leq N_0 $. Decompositions of $ \alpha_n^{(N)} $ and $ \widetilde{\alpha}_n^{(N)} $ respectively given by~\eqref{DecompalphN} and~\eqref{DecompalphtN} imply that
 \begin{align}
    &\widetilde{\alpha}_n^{(N)}(f) - \alpha_n^{(N)}(f) \notag \\
    &= \sum_{j=1}^{m_N} \frac{\mathbb{P}_{N}'(A_j^{(N)})}{\widetilde{\mathbb{P}}_n^{(N-1)}(A_j^{(N)})} (\widetilde{\alpha}_n^{(N-1)}(f \mathds{1}_{A_j^{(N)}}) - \alpha_n^{(N-1)}(f \mathds{1}_{A_j^{(N)}}))\notag \\
    &\quad + \alpha_n^{(N-1)}(f \mathds{1}_{A_j^{(N)}}) \left( \frac{\mathbb{P}_{N}'(A_j^{(N)})}{\widetilde{\mathbb{P}}_n^{(N-1)}(A_j^{(N)})}-\frac{P(A_j^{(N)})}{\mathbb{P}_n^{(N-1)}(A_j^{(N)})} \right) \notag \\
    &\quad - P(f \mathds{1}_{A_j^{(N)}}) \left(\frac{\widetilde{\alpha}_n^{(N-1)}(A_j^{(N)})}{\widetilde{\mathbb{P}}_n^{(N-1)}(A_j^{(N)})} - \frac{\alpha_n^{(N-1)}(A_j^{(N)})}{\mathbb{P}_n^{(N-1)}(A_j^{(N)})}\right) \notag \\
    &\quad+ \sqrt{\frac{n}{n_N}} \frac{P(f \mathds{1}_{A_j^{(N)}})}{\widetilde{\mathbb{P}}_n^{(N-1)}(A_j^{(N)})} \alpha_{N}'(A_j^{(N)}).\label{DiffTildN-N}
\end{align} By~\eqref{EqProcN=0} for $ N = 1 $ we have in particular \begin{align*}
    \widetilde{\alpha}_n^{(1)}(f) - \alpha_n^{(1)}(f)&= \sum_{j=1}^{m_1}  \alpha_n(f\mathds{1}_{A_j^{(1)}}) \left( \frac{\mathbb{P}_{n_1}'(A_j^{(1)})-P(A_j^{(1)})}{\mathbb{P}_n(A_j^{(1)})} \right)\\
    &\quad +\sqrt{\frac{n}{n_1}} \frac{P(f\mathds{1}_{A_j^{(1)}})}{\mathbb{P}_n(A_j^{(1)})} \alpha_{n_1}'(A_j^{(1)}),
\end{align*} which is uniformly and roughly bounded by \begin{align}
    ||\widetilde{\alpha}_n^{(1)}-\alpha_n^{(1)}||_\mathcal{F} \leq \frac{m_{(N)} K_\mathcal{F} \Lambda_n' }{p_{(N)} - \Lambda_n/\sqrt{n}} \sqrt{\frac{n}{n_{(N)}}} (1+\Lambda_n/\sqrt{n}).\label{BoundTild1-1}
\end{align} Let $ C_{n,N} =4m_{(N)} K_\mathcal{F}/(p_{(N)}-\Lambda_n/\sqrt{n})^2 $. Equality~\eqref{DiffTildN-N} implies also
 \begin{align*}
    &||\widetilde{\alpha}_n^{(N)}(f) - \alpha_n^{(N)}(f)||_\mathcal{F} \\
    &\leq C_{n,N} \left( ||\widetilde{\alpha}_n^{(N-1)} - \alpha_n^{(N-1)}||_\mathcal{F} + \frac{\Lambda_n^2}{\sqrt{n}} + \frac{\Lambda_n'(\Lambda_n + \sqrt{n})}{\sqrt{n_{(N)}}}\right).
\end{align*}
By induction of the last inequality and noticing that for all $ n >0$, $ m_{(N)} K_\mathcal{F}/(p_{(N)}-\Lambda_n/\sqrt{n})^2 \geq 1 $, we have
 \begin{align*}
    ||\widetilde{\alpha}_n^{(N)}(f) - \alpha_n^{(N)}(f)||_\mathcal{F} &\leq C_{n,N}^{N-1}  ||\widetilde{\alpha}_n^{(1)}-\alpha_n^{(1)}||_\mathcal{F} \\
    &\quad + (N-1)C_{n,N}^{N-1}\left(\frac{\Lambda_n^2}{\sqrt{n}} + \frac{\Lambda_n'(\Lambda_n+\sqrt{n}) }{\sqrt{n_{(N)}}} \right) ,
\end{align*}
then inequality~\eqref{BoundTild1-1} immediately implies that \begin{align*}
     ||\widetilde{\alpha}_n^{(N)}(f) - \alpha_n^{(N)}(f)||_\mathcal{F} \leq  N C_{n,N}^{N}\left( \frac{\Lambda_n^2}{\sqrt{n}} + \frac{\Lambda_n'(\Lambda_n +\sqrt{n})}{\sqrt{n_{(N)}}} \right).
 \end{align*}Since the right-hand side of the last inequality is increasing with $ N $ we find that for all $ t > 0 $,
  \begin{align}
    &\mathbb{P}\left(\sup_{1 \leq N \leq N_0} ||\widetilde{\alpha}_n^{(N)}(f) - \alpha_n^{(N)}(f)||_\mathcal{F} > t \right) \notag \\
    &\leq \mathbb{P}\left(\frac{C_0}{(p_{(N_0)}-\Lambda_n/\sqrt{n})^{2N_0}}  \left( \frac{\Lambda_n^2}{\sqrt{n}}+\frac{\Lambda_n'(\Lambda_n+\sqrt{n}) }{\sqrt{n_{(N_0)}}} \right) > t\right),\label{BoundProbN-N}
\end{align}
with $ C_0 = N_0^2(4m_{(N_0)} K_\mathcal{F}N_0)^{N_0}>0 $. There exists $ n_2 > 0 $ such that for all $ n > n_2 $ it holds $D\sqrt{\log(n)/n} \leq p_{(N_0)}/2 \leq 1/2 $. For $ n > n_2 $ we have according to~\eqref{BoundDevalphaNAndalphatN} and~\eqref{BoundProbN-N}, \begin{align*}
    &\mathbb{P}\left(\sup_{1 \leq N \leq N_0} ||\widetilde{\alpha}_n^{(N)}(f) - \alpha_n^{(N)}(f)||_\mathcal{F} > t \right) \\
    &\leq \mathbb{P}\left(\Lambda_n > D\sqrt{\log(n)}\right) \\
    &\quad +\mathbb{P}\left(\Lambda_n' > \frac{1}{2}\sqrt{\frac{n_{(N_0)}}{n}} \left(\frac{t p_{(N_0)}^{2N_0}}{4^{N_0} C_0}-\frac{D^2 \log(n)}{\sqrt{n}} \right) \right) \\
    &\leq \frac{1}{3n^2} +\mathbb{P}\left(\Lambda_n' > \frac{1}{2}\sqrt{\frac{n_{(N_0)}}{n}} \left(\frac{t p_{(N_0)}^{2N_0}}{4^{N_0} C_0}-\frac{D^2 \log(n)}{\sqrt{n}} \right) \right).
\end{align*} By using~\eqref{BoundDevalphaNAndalphatN} again, the last inequality implies $$ \mathbb{P}\left(\sup_{1 \leq N \leq N_0} ||\widetilde{\alpha}_n^{(N)}(f) - \alpha_n^{(N)}(f)||_\mathcal{F} > t_n \right) \leq \frac{2}{3n^2}, $$ for all $ n > n_2 $ and $$ t_n = \frac{4^{N_0+1}C_0 D}{p_{(N_0)}^{2N_0}}\left(\sqrt{\frac{n\log(n)}{n_{(N_0)}}}  + \frac{D\log(n)}{\sqrt{n}} \right). $$ By definition of $ v_n $, there exists $ d_2 > \max(d_1,4^{N_0+1}C_0 D/p_{(N_0)}^{2N_0}) $ and $ n_3 > 0 $ such that for all $ n > n_3 $, $$ d_2\left(v_n + \sqrt{\frac{n\log(n)}{n_{(N_0)}}}\right) > d_1 v_n + t_n. $$ Then~\eqref{ThmApproxForteResult} is proved for $ d_0=d_2 $ and $ n_0 = \max(n_0,n_1,n_3) $.

\subsection{Proof of Proposition~\ref{pro_Ztest}}

\noindent According to Theorem 2.1 of~\cite{AlbBert18} and Theorem~\ref{ThmApproxForte}, we can construct i.i.d random variables $ X_1,\dots,X_n $ with law $ P $ and $ z_n \sim \mathcal{N}(0,1) $ such that for $ n>n_1 $ for some $ n_1 > 0 $, $ \mathbb{P}(\mathcal{Z}_n) \leq 1/n^2 $ with \begin{align*}
    \mathcal{Z}_n^{(N)} &= \left\{|\alpha_n(f)/\sigma_f-z_n| > u_n \right\}\bigcup \left\{ |\alpha_n^{(N)}(f)/\sigma_f^{(N)}-z_n| > u_n \right\},
\end{align*} where $ u_n $ is a sequence with null limit. The strong approximation implies that \begin{align}
    \lim_{n \to +\infty} \frac{\mathbb{P}(|Z_n| \leq t_\alpha)}{\mathbb{P}(|z_n + M_n/\sigma_f| \leq t_\alpha)} = 1, \quad \lim_{n \to +\infty} \frac{\mathbb{P}(|Z_n^{(N)}| \leq t_\alpha)}{\mathbb{P}(|z_n + M_n/\sigma_f^{(N)}| \leq t_\alpha)} = 1, \label{Lim_pro_Ztest}
\end{align} with $ M_n = \sqrt{n}(P(f)-P_0(f)) $.If we denote $ f_{\mu,\sigma^2} $ the density function of $ \mathcal{N}(\mu,\sigma^2) $ then \begin{align*}
    \mathbb{P}(|z_n + M_n/\sigma_f| \leq t_\alpha) &\geq 2t_\alpha \inf_{[-t_\alpha,t_\alpha]} f_{M_n, 1} \\
    &\geq \frac{2t_\alpha}{\sqrt{2\pi}} \exp\left(-( M_n/\sigma_f+t_\alpha)^2 \right), \\
    \mathbb{P}(|z_n + M_n/\sigma_f^{(N)}| \leq t_\alpha) &\leq 2t_\alpha \sup_{[-t_\alpha,t_\alpha]} f_{M_n,1} \\
    &\leq \frac{2t_\alpha}{\sqrt{2\pi}} \exp \left(-(M_n/\sigma_f^{(N)}-t_\alpha)^2 \right).
\end{align*} which implies \begin{align*}
    \frac{\mathbb{P}(|z_n + M_n/\sigma_f| \leq t_\alpha)}{\mathbb{P}(|z_n + M_n/\sigma_f^{(N)}| \leq t_\alpha)} > \exp\left(M_n^2\left(\frac{1}{\sigma_f^{(N)}} - \frac{1}{\sigma_f} \right) -2t_\alpha |M_n|\left(\frac{1}{\sigma_f^{(N)}} + \frac{1}{\sigma_f} \right) \right)
\end{align*} For $ n $ large enough \begin{align}
    \frac{\mathbb{P}(|z_n + M_n/\sigma_f| \leq t_\alpha)}{\mathbb{P}(|z_n + M_n/\sigma_f^{(N)}| \leq t_\alpha)} \geq \exp\left(M_n^2\left(\frac{1}{\sigma_f^{(N)}} - \frac{1}{\sigma_f} \right)  \right). \label{Bound_pro_Ztest}
\end{align} Then~\eqref{Lim_pro_Ztest} and~\eqref{Bound_pro_Ztest} imply~\eqref{pro_Ztest_3}.

\subsection{Proof of Proposition~\ref{pro_chi2}}

Denote $ X\cdot Y $ the product scalar of $ X $ and $ Y $ and $ \mathcal{C} \in \mathbb{R}^m $ the random vector defined by $$ \mathcal{C} = (C_1,\dots,C_m)= (\mathds{1}_{B_1}/\sqrt{P(B_1)},\dots,\mathds{1}_{B_m}/\sqrt{P(B_m)}) .$$ We deal with the case $ (H_0) $ at Step 1 and the case $ (H_1) $ at Step 2.\medskip

\noindent \textbf{Step 1.} Under $ (H_0) $, $ T_n =  \alpha_n[\mathcal{C}] \cdot \alpha_n[\mathcal{C}]^T $, $ T_n^{(N)}  = \alpha_n^{(N)}[\mathcal{C}] \cdot \alpha_n^{(N)}[\mathcal{C}]^T $ and $ \widetilde{T}_n^{(N)} = \widetilde{\alpha}_n^{(N)}[\mathcal{C}] \cdot \widetilde{\alpha}_n^{(N)}[\mathcal{C}]^T $. Statistic $ T_n $ converges weakly to a multinormal random variable $ Y\sim \mathcal{N}(\mathbf{0},\Sigma) $ while $ T_n^{(N)},\widetilde{T}_n^{(N)} $ converge weakly to $ Y^{(N)} \sim \mathcal{N}(\mathbf{0},\Sigma^{(N)}) $ according to Theorem 2.1 of~\cite{AlbBert18} and Theorem~\ref{ThmApproxForte}. By Proposition 7 of~\cite{AlbBert18}, $ \Sigma-\Sigma^{(N)} $ is positive definite which implies for all $ \alpha > 0 $, $$ \mathbb{P}(Y\cdot Y^T\geq t_\alpha) \geq \mathbb{P}(Y^{(N)}\cdot (Y^{(N)})^T\geq t_\alpha), $$ and consequently~\eqref{pro_chi2_1},~\eqref{pro_chi2_2} by definition of weak convergence.\medskip

\noindent \textbf{Step 2.} Under $ (H_1) $, there exists $ i \in \{1,\dots,m\} $ such that $ P_0(B_i)\neq P(B_i) $ which implies \begin{align*}
    \min(|T_n|,|T_n^{(N)}|,\widetilde{T}_n^{(N)}|) &> -\Lambda_n^2 -2\sqrt{n} \Lambda_n |P_0(C_i)-P(C_i)| \\
    &\quad+ n(P_0(C_i)-P(C_i))^2.
\end{align*} By Borel-Cantelli and~\eqref{BoundDevalphaNAndalphatN} with probability one there exists $ n_1 > 0 $ such that for all $ n > n_1 $, $ \Lambda_n < D\sqrt{\log(n)} $. For $ n > n_1 $, we have \begin{align*}
    t_n &< \min(|T_n|,|T_n^{(N)}|,\widetilde{T}_n^{(N)}|),  \\
    t_n &= -D^2 \log(n)-2D \sqrt{n\log(n)}|P_0(C_i)-P(C_i)| \\
    &\quad+ n(P_0(C_i)-P(C_i))^2.
\end{align*}  Since $ \lim_{n \to +\infty} t_n = +\infty $, for all $ \alpha \in (0,1) $ there exists $ n_2 > 0 $ such that $ t_n > t_\alpha $ for all $ n > n_2 $. Inequality~\eqref{pro_chi2_3} is satisfied for $ n_0 = \max(n_1,n_2) $.

\appendix

\section{Numerical example of a raked mean}
\label{Appendix0}

The usual way to calculate the mean of $ X_1,\dots,X_n $ is to sum the data $ X_i $ multiplied by the weights $ w_i= 1/n $. If we have the auxiliary information $ P[\mathcal{A}^{(N)}] = (P(A_1^{(N)}),\dots,P(A_{m_N}^{(N)})) $ for $ 1 \leq N \leq N_0 $ we want to change iteratively the initial weights $w_i $ in new weights $ w_i^{(N)} $ such that $ \sum_{i=1}^n w_i^{(N)} $ and $$ \sum_{i=1}^n w_i^{(N)} \mathds{1}_{A_j}^{(N)}(X_i) = P(A_j^{(N)}), $$ for any $ 1 \leq N \leq N_0 $ and $ 1 \leq j \leq m_N $. Recall that it does not imply that $ \sum_{i=1}^n w_i^{(N_1)} \mathds{1}_{A_j}^{(N_2)}(X_i) = P(A_j^{(N_2)}) $ with $ N_1\neq N_2 $ and $ 1 \leq j \leq N_2 $. For this example one takes $ N_0=2, \mathcal{A}^{(2)} = \{A_1,A_2,A_3\}, \mathcal{B}=\{B_1,B_2\} $ and one generates normal random values $ X_i $ with fixed variances $ \sigma^2=0.1 $ and such that the probabilities and conditional expectations are given by the following table: \begin{center}
    \begin{tabular}{|c|c|c|c|}
     \hline $ P(A_i \cap B_j) $ & $ A_1 $ & $ A_2 $ & $ A_3 $  \\ \hline
     $ B_1 $ & 0.2 & 0.25 & 0.1   \\ \hline
     $ B_2 $ & 0.25 &0.1 & 0.1 \\  \hline
\end{tabular}
    \captionof{table}{Probabilities of sets}
\end{center}
\begin{center}
    \begin{tabular}{|c|c|c|c|}
     \hline $ \mathbb{E}[X|A_i \cap B_j] $ & $ A_1 $ & $ A_2 $ & $ A_3 $  \\ \hline
     $ B_1 $ & 0.75& -0.5 & 1  \\ \hline
     $ B_2 $ &0.5 &0.25 &-0.5 \\  \hline
\end{tabular}
    \captionof{table}{Conditional expectations of the generated random variables}
\end{center} In particular, \begin{align*}
    P[\mathcal{A}] &= (P(A_1),P(A_2),P(A_3)) = (0.45,0.35,0.2), \\
    P[\mathcal{B}] &= (P(B_1),P(B_2)) = (0.55,0.45), \\
    P(X) &=\mathbb{E}[X] = 0.225, \\
    \mathbb{E}[X|\mathcal{A}] &= (\mathbb{E}[X|A_1],\mathbb{E}[X|A_2],\mathbb{E}[X|A_3]) \simeq (0.611, -0.286, 0.25), \\
    \mathbb{E}[X|\mathcal{B}] &= (\mathbb{E}[X|B_1],\mathbb{E}[X|B_2]) \simeq (0.227, 0.222).
\end{align*} We generate $ n=10 $ values and we obtain the following data \begin{center}
    \begin{tabular}{|c|c|c|} \hline
         $X_i$ & $ \mathcal{A} $ & $ \mathcal{B} $  \\ \hline
         0.953 & 1&1 \\ \hline
         0.975 & 1 & 1 \\ \hline
         0.058 & 1&1 \\ \hline
         -0.766 & 2&1 \\ \hline
         -0.644 & 2 &1\\ \hline
         -0.819 & 2&1 \\ \hline
         0.028 & 2 &2\\ \hline
         0.627 & 2 &2\\ \hline
         1.04 & 3 &1\\ \hline
         -0.904 & 3 & 2 \\ \hline
    \end{tabular}
    \captionof{table}{Generated random variables}
\end{center} In this case, the usual mean is the sum of all $ X_i $ over 10 that is we assign the weight $1/n=0.1 $ at each $ X_i $ and we have $ \mathbb{P}_n(X) \simeq 0.055 $. When we rake one time we assign the weights $ 0.15,0.07,0.1 $ at individuals belonging respectively to $ A_1,A_2,A_3 $. The raked mean for $ N=1 $ is \begin{align*}
    \mathbb{P}_n^{(1)}(X) &= 0.15 \times \frac{P(A_1)}{\mathbb{P}_n(A_1)} + 0.07 \times \frac{P(A_2)}{\mathbb{P}_n(A_2)} + 0.1 \times \frac{P(A_3)}{\mathbb{P}_n(A_3)} \simeq 0.2.
\end{align*} When the algorithm is stabilized in this case the final weights are given by the following table: \begin{center}
     \begin{tabular}{|c|c|c|c|}
     \hline $ w_i^{(\infty)} $ & $ A_1 $ & $ A_2 $ & $ A_3 $  \\ \hline
     $ B_1 $ & 0.15 & 0.024 & 0.029   \\ \hline
     $ B_2 $ & X &0.139 & 0.17 \\  \hline
\end{tabular}
    \captionof{table}{Final raked weights}
\end{center}
Notice that the cross means that we do not generate random variables belonging to $ A_1 \cap B_2 $ due to a low value of $ n $. The final raked mean is $ \mathbb{P}_n^{(\infty)}(X) \simeq 0.212 $ which is closer of $ P(X) $ than the usual mean $ \mathbb{P}_n(X) $.

\section{Calculation of $ \sigma_f^{(\infty)} $}\label{Appendix1}

We use the notations of the section 4.4 of~\cite{AlbBert18} concerning the proof of their Proposition 11 in the aim to establish the expression of $ \mathbb{G}^{(\infty)} $. The calculation uses the two following stochastic matrices
\begin{align}
\mathbf{P}_{\mathcal{A}|\mathcal{B}} &= \left( \begin{matrix}
        P(A|B) & P(A^C|B) \\
        P(A|B^C) & P(A^C|B^C)
    \end{matrix}\right)  = \left( \begin{matrix}
        p_{AB}/p_B & 1-p_{AB}/p_B \\
        (p_A-p_{AB})/p_{\overline{B}} & 1-(p_A-p_{AB})/p_{\overline{B}}
    \end{matrix}\right), \notag \\ \notag \\
    \mathbf{P}_{\mathcal{B}|\mathcal{A}} &= \left( \begin{matrix}
        P(B|A) & P(B^C|A) \\
        P(B|A^C) & P(B^C|A^C)
    \end{matrix}\right) = \left( \begin{matrix}
        p_{AB}/p_A & 1-p_{AB}/p_A \\
        (p_B-p_{AB})/p_{\overline{A}} & 1-(p_B-p_{AB})/p_{\overline{A}}
    \end{matrix}\right), \label{DefPABPBA}
\end{align} the two following conditional expectation vectors \begin{align}
    \mathbb{E}[f|\mathcal{A}] = (\mathbb{E}[f|A],\mathbb{E}[f|A^C]) , \quad \mathbb{E}[f|\mathcal{B}] = (\mathbb{E}[f|B],\mathbb{E}[f|B^C]),\label{DefEAEB}
\end{align} the two following covariance matrices \begin{align}
    \mathrm{Var}(\mathbb{G}[A]) = p_A p_{\overline{A}} \left( \begin{matrix}
        1&-1 \\ -1&1
    \end{matrix} \right),  \quad \mathrm{Var}(\mathbb{G}[B]) = p_B p_{\overline{B}} \left( \begin{matrix}
        1&-1 \\ -1&1
    \end{matrix} \right). \label{DefCovGAGB}
\end{align} and the two following vectors \begin{align*}
    V_1(f) &= \mathbb{E}[f|\mathcal{A}]-\mathbf{P}_{\mathcal{B}|\mathcal{A}} \cdot \mathbb{E}[f|\mathcal{B}] \\
    &= \left( \begin{matrix}
        \mathbb{E}[f|A] \\
        \mathbb{E}[f|A^C]
    \end{matrix}\right) -  \left( \begin{matrix}
        p_{AB}/p_A & 1-p_{AB}/p_A \\
        (p_B-p_{AB})/p_{\overline{A}} & 1-(p_B-p_{AB})/p_{\overline{A}}
    \end{matrix}\right) \cdot
    \left( \begin{matrix}
        \mathbb{E}[f|B] \\
        \mathbb{E}[f|B^C]
    \end{matrix}\right) \\
    &= (\mathbb{E}[f](p_A-p_{AB})-\mathbb{E}[f|A]p_A p_{\overline{B}} +\mathbb{E}[f|B](p_{AB}-p_Ap_B))\cdot \left( \begin{matrix}
        -1/p_A p_{\overline{B}} \\
        1/p_{\overline{A}} p_{\overline{B}}
    \end{matrix}\right) , \\ \\
    V_2(f) &= \mathbb{E}[f|\mathcal{B}]-\mathbf{P}_{\mathcal{A}|\mathcal{B}} \cdot \mathbb{E}[f|\mathcal{A}] \\
    &=\left( \begin{matrix}
        \mathbb{E}[f|B] \\
        \mathbb{E}[f|B^C]
    \end{matrix}\right) -  \left( \begin{matrix}
        p_{AB}/p_B & 1-p_{AB}/p_B \\
        (p_A-p_{AB})/p_{\overline{B}} & 1-(p_A-p_{AB})/p_{\overline{B}}
    \end{matrix}\right) \cdot
    \left( \begin{matrix}
        \mathbb{E}[f|A] \\
        \mathbb{E}[f|A^C]
    \end{matrix}\right) \\
    &= (\mathbb{E}[f](p_B-p_{AB})-\mathbb{E}[f|B]p_{\overline{A}} p_{B} +\mathbb{E}[f|A](p_{AB}-p_Ap_B))\cdot \left( \begin{matrix}
        -1/p_{\overline{A}} p_B \\
        1/p_{\overline{A}} p_{\overline{B}}
    \end{matrix}\right).
\end{align*} The eigenvalues of $ \mathbf{P}_{\mathcal{A}|\mathcal{B}} \cdot \mathbf{P}_{\mathcal{B}|\mathcal{A}} $ and $ \mathbf{P}_{\mathcal{B}|\mathcal{A}} \cdot \mathbf{P}_{\mathcal{A}|\mathcal{B}} $ are 1 and $ T_1=T_2 = (p_{AB}-p_A p_B)^2/p_A p_{\overline{A}} p_B p_{\overline{B}} $. Their eigenvectors associated to $ T_1 $ and $ T_2 $ are respectively  $ (p_{\overline{B}}/p_B, -1)^t $ and $ (p_{\overline{A}}/p_A, -1)^t $ which implies $$ U_1 = \left( \begin{matrix}
    1& p_{\overline{A}}/p_A \\
    1 &-1
    \end{matrix} \right), \quad U_2 = \left( \begin{matrix}
    1& p_{\overline{B}}/p_B \\
    1 &-1
    \end{matrix} \right). $$ For the case of two marginals, Albertus and Berthet showed that $ \mathbb{G}^{(N)} $ converge almost surely to $ \mathbb{G}^{(\infty)}(f) = \mathbb{G}(f) - S_{1,even}(f)^t \cdot \mathbb{G}[\mathcal{A}] - S_{2,odd}(f)^t \cdot \mathbb{G}[\mathcal{B}] $ where
\begin{align*}
    S_{1,even}(f) &= U_1 \left( \begin{matrix}
        0 & 0 \\
        0 & (1-T_1)^{-1}
    \end{matrix} \right) \cdot U_1^{-1} \cdot  V_1(f) = C_{1,even}(f) \left( \begin{matrix}
        -p_{\overline{A}}p_B \\
        p_Ap_B
    \end{matrix}\right), \\
    C_{1,even}(f) &= \frac{\mathbb{E}[f|B](p_{AB}-p_A p_B)-\mathbb{E}[f|A]p_Ap_{\overline{B}} - \mathbb{E}[f](p_{AB}-p_A)}{p_Ap_B p_{\overline{A}}p_{\overline{B}} - (p_{AB}-p_A p_B)^2}, \\ \\
    S_{2,odd}(f) &= U_2 \left( \begin{matrix}
        0 & 0 \\
        0 & (1-T_2)^{-1}
    \end{matrix} \right) \cdot U_2^{-1} \cdot  V_2(f) = C_{2,odd}(f) \left( \begin{matrix}
        -p_{A}p_{\overline{B}} \\
        p_Ap_B
    \end{matrix} \right), \\
    C_{2,odd}(f) &= \frac{ \mathbb{E}[f|A](p_{AB}-p_A p_B)-\mathbb{E}[f|B]p_{\overline{A}} p_B- \mathbb{E}[f](p_{AB}-p_B)}{p_Ap_B p_{\overline{A}}p_{\overline{B}} - (p_{AB}-p_A p_B)^2}.
\end{align*}  By linearity of $ f \mapsto \mathbb{G}(f) $ and the fact that $ \mathbb{G}(a) = 0 $ for any constant $ a \in \mathbb{R} $ one can write $$ \mathbb{G}^{(\infty)}(f) = \mathbb{G}\left(f + p_B C_{1,even}(f) \mathds{1}_A + p_A C_{2,odd}(f) \mathds{1}_B \right), $$ which implies that \begin{align*}
    \sigma_f^{(\infty)} &= \mathrm{Var}(\mathbb{G}^{(\infty)}(f)) \\
    &= \mathrm{Var}(f) + \mathrm{Var}(p_B C_{1,even}(f) \mathds{1}_A + p_A C_{2,odd}(f) \mathds{1}_B) \\
    &\qquad +2\mathrm{Cov}(f, p_B C_{1,even}(f) \mathds{1}_A + p_A C_{2,odd}(f) \mathds{1}_B) \\
    &= \mathrm{Var}(f) + p_A p_{\overline{A}} p_B^2 C_{1,even}^2(f) + p_A^2 p_B p_{\overline{B}} C_{2,odd}^2(f) \\
    &\qquad + 2p_A p_B C_{1,even}(f) C_{2,odd}(f)(p_{AB}-p_A p_B) \\
    &\qquad +2p_A p_B\left( C_{1,even}(f) \Delta_A + C_{2,odd}(f) \Delta_B \right)
\end{align*} With some calculations we find the simple expression of $ \sigma_f^{(\infty)} $ given by~\eqref{DefSigInf}.

\end{document}